\documentclass[12pt]{amsproc}
\usepackage[cp1251]{inputenc}
\usepackage[T2A]{fontenc} 
\usepackage[russian]{babel} 
\usepackage{amsbsy,amsmath,amssymb,amsthm,graphicx,color} 
\usepackage{amscd}
\def\le{\leqslant}
\def\ge{\geqslant}

\newtheorem{prop}{Предложение}

\theoremstyle{definition}

\theoremstyle{remark}

\begin {document}
\unitlength=1mm
\title[Формулы Морделла и Туэнтера]
{Формулы Морделла и Туэнтера, функции Куберта и теоремы взаимности для сумм Апостола-Дедекинда}
\author{Г. Г. Ильюта}
\email{ilgena@rambler.ru}
\address{}
\thanks{Работа поддержана грантом РФФИ-20-01-00579}

\bigskip

\begin{abstract}
Мы свяжем суммы Дедекинда и некоторые формулы для числовых полугрупп.

We connect Dedekind sums and some formulas for numerical semigroups.
\end{abstract}

\maketitle
\tableofcontents

\bigskip

\section{Введение}

\bigskip

  Теоремы взаимности являются ключевыми фактами для понимания структуры различных вариаций сумм Дедекинда и близких к ним сумм, в частности, они существенно ускоряют вычисление этих сумм, поскольку позволяют запустить алгоритм Евклида. Например, для сумм Апостола-Дедекинда
$$
s_n(a,b):=\sum_{k=1}^{b-1}\frac{k}{b}\bar B_n\left(\frac{ka}{b}\right)   
$$
теорема взаимности (для нечётного $n$) имеет вид \cite{1}
$$
ab^ns_n(a,b)+ba^ns_n(b,a)
$$
$$
=\sum_{k=0}^{n+1}\frac{(-1)^kn!}{k!(n-k+1)!}B_kB_{n+1-k}a^kb^{n-k+1}+\frac{nB_{n+1}}{n+1}                                              \eqno (1)
$$
Здесь $a$ и $b$ -- взаимно простые натуральные числа, а функции Бернулли $\bar B_n(x)$ определяются равенством
$$
\bar B_n(x):=B_n(x-[x]),
$$
где $[x]$ -- целая чать $x$,
$$
B_n(x)=\sum_{i=0}^n\binom{n}{i}B_ix^{n-i}        \eqno (2)
$$
-- многочлены Бернулли, $B_n=B_n(0)$ -- числа Бернулли. Ниже мы без дополнительного упоминания используем равенства $\bar B_n(x+1)=\bar B_n(x)$ для всех $x$ и $\bar B_n(x)=B_n(x)$ для $0\le x<1$.

  Функции Куберта \cite{10}, \cite{13} (известные также как репликативные функции \cite{12}, 1.2.4, Ex. 39, 40) по определению удовлетворяют соотношению вида
$$
f(ax)=\alpha_a\sum_{k=1}^{a-1}f\left(x+\frac{k}{a}\right)+\beta_a.    \eqno (3)
$$
В \cite{1} соотношение вида (3) для функций Бернулли
$$
\bar B_n(ax)=a^{n-1}\sum_{k=1}^{a-1}\bar B_n\left(x+\frac{k}{a}\right)         
$$
позволяет при доказательстве формулы (1) перейти от суммирования по одному индексу в определении сумм Апостола-Дедекинда к симметричному по $a$ и $b$ суммированию по двум индексам. В п.3 мы интерпретируем этот приём как переход от суммирования по множеству Апери $Ap_b(S_{a,b})$ числовой полугруппы $S_{a,b}$, порождённой числами $a$ и $b$, к суммированию по множеству Апери $Ap_{ab}(S_{a,b})$, что позволит применить формулу Морделла (8) (см. в следующем абзаце объяснение этой терминологии). В частности, в Предложении 3 мы получим другое представление для правой части формулы (1), которое содержит суммы Сильвестра числовой полугруппы $S_{a,b}$. Также в Предложении 3 представлена общая схема использования формулы Морделла при выводе теорем взаимности для сумм вида $\sum(k/b)f(ka/b)$, где функция $f(x)$ берётся из некоторого класса функций, содержащего функции Куберта (соотношение (3) нам понадобится не для всей области определения переменной $x$, а только для нескольких рациональных значений этой переменной, см. условия Предложения 3).

  Множество $S\subset \mathbb Z_{\ge 0}$ называется числовой полугруппой, если оно содержит $0$, замкнуто относительно сложения и имеет конечное дополнение в $\mathbb Z_{\ge 0}$. Для $m\in S$, $m\neq 0$, множество Апери $Ap_m(S)$ определяется следующим образом:
$$
Ap_m(S):=\{s\in S:s-m\notin S\}.
$$
Другими словами, $Ap_m(S)$ состоит из минимальных элементов в пересечениях полугруппы $S$ с классами вычетов по $\mod m$. Нам понадобится следующая формула Гассерта-Шора \cite{8}: для произвольной функции $f$ на $\mathbb Z_{\ge 0}$
$$
\sum_{k\in\mathbb Z_{\ge 0}\setminus S}(f(k+m)-f(k))=\sum_{k\in Ap_m(S)}f(k)-\sum_{k=0}^{m-1}f(k).         \eqno (4)
$$
Как отмечалось в \cite{8}, частным случаем формулы Гассерта-Шора (4), отвечающим множеству $Ap_b(S_{a,b})$, является формула Туэнтера \cite{16}
$$
\sum_{k\in C_{a,b}}(f(k+b)-f(k))=\sum_{k=0}^{b-1}(f(ak)-f(k)),         \eqno (5)
$$
где $C_{a,b}:=\mathbb Z_{\ge 0}\setminus S_{a,b}$. Это следует из равенства       
$$
Ap_b(S_{a,b})=\{0a,1a,\dots,(b-1)a\}.        \eqno (6)
$$   
Мы покажем, что частный случай формулы Гассерта-Шора (4), отвечающий множеству $Ap_{ab}(S_{a,b})$, фактически содержится (совпадает с разностью формул (7) и (8), мы приводим их в п.2) в статье Л. Морделла, посвящённой формуле взаимности для сумм Дедекинда и вышедшей в 1951 году \cite{14}. Это следует из равенства 
$$
Ap_{ab}(S_{a,b})=\{ia+jb:i=0,1,\dots,b-1;j=0,1,\dots,a-1\},       \eqno (7)
$$
которое мы докажем в п.2. Назовём этот частный случай формулы Гассерта-Шора (4) формулой Морделла
$$
\sum_{k\in C_{a,b}}(f(k+ab)-f(k))=\sum_{i=0}^{b-1}\sum_{j=0}^{a-1}f(ia+jb)-\sum_{k=0}^{ab-1}f(k).         \eqno (8)
$$

  Обозначим через $R_n(a,b)$ множество всех элементов числовой полугруппы $S_{a,b}$, имеющих ровно $n$ представлений в виде линейной комбинации образующих. В Предложении 2 мы покажем как множество $Ap_{ab}(S_{a,b})$ определяет множества $R_n(a,b)$. Этот факт получен при сравнении описаний множества $Ap_{ab}(S_{a,b})$ в Предложении 1 и множеств $R_n(a,b)$ в \cite{2}. В частности, множество $Ap_{ab}(S_{a,b})$ состоит из всех элементов числовой полугруппы $S_{a,b}$, имеющих единственное представление в виде линейной комбинации образующих
$$
Ap_{ab}(S_{a,b})=R_1(a,b).
$$

  Мы рассматриваем формулу Морделла и формулу Туэнтера как пути к теоремам взаимности для некоторых классов сумм. О формуле Морделла мы уже говорили (её применение обсуждается в п.3), а в п.4 аналогичным образом рассматривается формула Туэнтера. В Предложении 4 представлена общая схема её применения и в качестве примера получено ещё одно (также содержащее суммы Сильвестра) представление для правой части формулы (1).

  Возможно, что правильно было бы рассматривать каждый элемент числовой полугруппы $S_{a,b}$ и отвечающий ему частный случай формулы Гассерта-Шора (4) как путь к теоремам взаимности для некоторого класса сумм, но мы не будем заниматься здесь этой темой.

  Интересным инвариантом изолированной особой точки голоморфной функции является её спектр -- конечный набор рациональных чисел. Например, числа $e^{2\pi ik}$, $k$ пробегает спектр особенности, образуют спектр оператора монодромии особенности. В \cite{4}, \cite{5}, \cite{7} изучались (сдвинутые) степенные суммы спектральных чисел. В п.5 мы получим две формулы для (сдвинутых) степенных сумм
$$
\sum_{i=1}^{b-1}\sum_{j=1}^{a-1}\left(\frac{x}{ab}+\frac{j}{a}+\frac{i}{b}\right)^n
$$
спектральных чисел особенности Брискорна $z_1^a+z_2^b$.

  Для узла в сфере $S^3$, его матрицы Зейферта $Q$ и числа $z\in\mathbb C$, $|z|=1$, $z\ne 1$, обозначим через $\sigma(z)$ сигнатуру эрмитовой формы
$$
(1-z)Q+(1-\bar z)Q^t.
$$
Известно, что функция $\sigma(z)$ не зависит от выбора матрицы Зейферта узла. Для $(a,b)$-торического узла в \cite{11} доказана формула
$$
\int_0^1\sigma(e^{2\pi ix})dx=-\frac{1}{3}(a-\frac{1}{a})(b-\frac{1}{b}).  \eqno (9)
$$
В \cite{3} получено элементарное доказательство этого равенства, которое в п.6 мы сделаем более коротким. 

  Наши доказательства, относящиеся к двум предыдущим абзацам (п.5 и п.6), используют формулу (10) (формулу (7) из \cite{14}), в которой присутствует суммирование по множеству $Ap_{ab}(S_{a,b})$. Связь сигнатуры торического узла и классических сумм Дедекинда (сумм по множеству $Ap_b(S_{a,b})$) известна уже более полувека (см. список литературы в \cite{3}). Торический узел является пересечением особого слоя функции $z_1^a+z_2^b$ и сферы с центром в особой точке в $\mathbb C^2$. Возможно, что в этом круге тем, как и в случае теорем взаимности, окажутся полезными суммы по множеству $Ap_m(S_{a,b})$ для произвольного элемента $m$ числовой полугруппы $S_{a,b}$. Формула Гассерта-Шора сводит все такие суммы к суммам по дополнению $C_{a,b}$ к числовой полугруппе $S_{a,b}$ в $\mathbb Z_{\ge0}$.

\bigskip

\section{Множества Апери и формула Гассерта-Шора}

\bigskip

  В следующем Предложении мы покажем, что, зная множество Апери одного элемента $m_1\neq 0$ числoвой полугруппы $S$, можно легко найти множество Апери любого другого элемента $m_2\neq 0$.

\begin{prop}\label{prop1} Для любых не равных $0$ элементов $m_1,m_2$ числoвой полугруппы $S$
$$
\frac{\sum_{k\in Ap_{m_1}(S)}q^k}{q^{m_1}-1}=\frac{\sum_{k\in Ap_{m_2}(S)}q^k}{q^{m_2}-1}
$$
$$
=\sum_{k\in\mathbb Z_{\ge 0}\setminus S}q^k+\frac{1}{q-1}.
$$
В частности,
$$
Ap_{ab}(S_{a,b})=\{ia+jb:i=0,1,\dots,b-1;j=0,1,\dots,a-1\}.
$$
\end{prop}

  Доказательство. Сравниваем формулу Гассерта-Шора для $m_1$ и $f(k)=q^k$
$$
(q^{m_1}-1)\sum_{k\in\mathbb Z_{\ge 0}\setminus S}q^k=\sum_{k\in Ap_{m_1}(S)}q^k-\frac{q^{m_1}-1}{q-1}
$$
с аналогичной формулой для $m_2$. Полагая $m_1=ab$, $m_2=b$, получим
$$
\sum_{k\in Ap_{ab}(S_{a,b})}q^k=(q^{ab}-1)\frac{\sum_{k=0}^{b-1}q^{ak}}{q^b-1}
$$
$$
=\frac{(q^{ab}-1)^2}{(q^a-1)(q^b-1)}
=\sum_{k\in\{ia+jb:i=0,1,\dots,b-1;j=0,1,\dots,a-1\}}q^k.
$$

\begin{prop}\label{prop2} Для $n\ge 1$
$$
\{k+(n-1)ab:k\in Ap_{ab}(S_{a,b})\}=R_n(a,b),
$$
\end{prop}

Доказательство. Сравниваем равенство (7) с  формулой из \cite{2} 
$$
\sum_{k\in R_n(a,b)}q^k=
\frac{q^{(n-1)ab}(q^{ab}-1)^2}{(q^a-1)(q^b-1)}.
$$

  Приведём формулы (7) и (8) из \cite{14}, разность которых совпадает с формулой Морделла (8):
$$
\sum_{i=1}^{b-1}\sum_{j=1}^{a-1}f(ia+jb)
$$
$$
=\sum_{\substack{0<i<b,0<j<a\\ia+jb<ab}}f(ia+jb)+\sum_{\substack{0<i<b,0<j<a\\ia+jb<ab}}f(2ab-ia-jb),    \eqno (10)
$$
$$
\sum_{k=1}^{ab-1}f(k)-\sum_{i=1}^{b-1}f(ia)-\sum_{j=1}^{a-1}f(jb)
$$
$$
=\sum_{\substack{0<i<b,0<j<a\\ia+jb<ab}}f(ia+jb)+\sum_{\substack{0<i<b,0<j<a\\ia+jb<ab}}f(ab-ia-jb).
$$
В \cite{14} доказано равенство
$$
C_{a,b}=\{ab-ia-jb:0<i<b,0<j<a,ia+jb<ab\}.
$$
Напомним результат Фробениуса $\max C_{a,b}=ab-a-b$. Формула Гассерта-Шора (4) очень легко сводится к определению числовой полугруппы и её множества Апери \cite{9}. Поэтому формула Туэнтера (5) практически эквивалентна равенству (6), а формула Морделла (8) -- равенству (7).

\bigskip

\section{Формула Морделла и теоремы взаимности}

\bigskip

  Определим суммы Сильвестра числовой полугруппы $S_{a,b}$ равенством
$$
S_{a,b}(m):=\sum_{k\in C_{a,b}}k^m.
$$
Для этих сумм известна формула \cite{15}
$$
S_{a,b}(m-1)=\sum_{i=0}^m\sum_{j=0}^{m-i}\frac{(m-1)!}{i!j!(m-i-j+1)!}B_iB_ja^{m-j}b^{m-i}-\frac{B_m}{m}.
$$
\begin{prop}\label{prop2} Пусть функция $f$ удовлетворяет соотношениям
$$
f\left(\frac{ia}{b}\right):=\alpha_a\sum_{j=0}^{a-1}f\left(\frac{i}{b}+\frac{j}{a}\right)+\beta_a,
$$
$$
f\left(\frac{jb}{a}\right):=\alpha_b\sum_{i=0}^{b-1}f\left(\frac{i}{b}+\frac{j}{a}\right)+\beta_b.
$$
Тогда
$$
\alpha_b\sum_{i=0}^{b-1}\frac{i}{b}f\left(\frac{ia}{b}\right)+\alpha_a\sum_{j=0}^{a-1}\frac{j}{a}f\left(\frac{jb}{a}\right)
$$
$$
=\alpha_a\alpha_b\sum_{k\in C_{a,b}}\left(\left(\frac{k}{ab}+1\right)f\left(\frac{k}{ab}+1\right)-\frac{k}{ab}f\left(\frac{k}{ab}\right)\right)
$$
$$
+\alpha_a\alpha_b\sum_{k=0}^{ab-1}\frac{k}{ab}f\left(\frac{k}{ab}\right)+\frac{\alpha_b\beta_a(b-1)}{2}+\frac{\alpha_a\beta_b(a-1)}{2},
$$
в частности, для $1$-периодической функции $f$ имеем
$$
\alpha_b\sum_{i=0}^{b-1}\frac{i}{b}f\left(\frac{ia}{b}\right)+\alpha_a\sum_{j=0}^{a-1}\frac{j}{a}f\left(\frac{jb}{a}\right)=\alpha_a\alpha_b\sum_{k\in C_{a,b}}f\left(\frac{k}{ab}\right)
$$
$$
+\alpha_a\alpha_b\sum_{k=0}^{ab-1}\frac{k}{ab}f\left(\frac{k}{ab}\right)+\frac{\alpha_b\beta_a(b-1)}{2}+\frac{\alpha_a\beta_b(a-1)}{2},
$$
для $f(x)=\bar B_n(x)$ имеем
$$
ab^ns_n(a,b)+ba^ns_n(b,a)=
$$
$$
=\sum_{k\in C_{a,b}}\sum_{i=0}^n\binom{n}{i}(ab)^iB_iS_{a,b}(n-i)+\sum_{i=0}^{n+1}\frac{n!(ab)^{n-i+1}B_iB_{n-i+1}}{i!(n-i+1)!}+\frac{nB_{n+1}}{n+1}.
$$
\end{prop}

  Доказательство. Применяя формулу Морделла (8), получим
$$
\alpha_b\sum_{i=0}^{b-1}\frac{i}{b}f\left(\frac{ia}{b}\right)+\alpha_a\sum_{j=0}^{a-1}\frac{j}{a}f\left(\frac{jb}{a}\right)
$$
$$
=\alpha_a\alpha_b\sum_{i=0}^{b-1}\sum_{j=0}^{a-1}\left(\frac{j}{a}+\frac{i}{b}\right)f\left(\frac{j}{a}+\frac{i}{b}\right)+\alpha_b\beta_a\sum_{i=0}^{b-1}\frac{i}{b}+\alpha_a\beta_b\sum_{j=0}^{a-1}\frac{j}{a}
$$
$$
=\alpha_a\alpha_b\sum_{k\in C_{a,b}}\left(\left(\frac{k}{ab}+1\right)f\left(\frac{k}{ab}+1\right)-\frac{k}{ab}f\left(\frac{k}{ab}\right)\right)
$$
$$
+\alpha_a\alpha_b\sum_{k=0}^{ab-1}\frac{k}{ab}f\left(\frac{k}{ab}\right)+\frac{\alpha_b\beta_a(b-1)}{2}+\frac{\alpha_a\beta_b(a-1)}{2}.
$$
Согласно \cite{1}, Lemma 1,
$$
m^n\sum_{k=1}^{m-1}\frac{k}{m}B_n\left(\frac{k}{m}\right)=\sum_{i=0}^{n+1}\frac{n!m^{n-i+1}B_iB_{n-i+1}}{i!(n-i+1)!}+\frac{nB_{n+1}}{n+1}              \eqno (11)
$$
Поэтому, используя равенство (2), получим
$$
ab^ns_n(a,b)+ba^ns_n(b,a)=(ab)^n\left(\sum_{k\in C_{a,b}}B_n\left(\frac{k}{ab}\right)+\sum_{k=0}^{ab-1}\frac{k}{ab}B_n\left(\frac{k}{ab}\right)\right)
$$
$$
=\sum_{k\in C_{a,b}}\sum_{i=0}^n\binom{n}{i}(ab)^iB_ik^{n-i}+\sum_{i=0}^{n+1}\frac{n!(ab)^{n-i+1}B_iB_{n-i+1}}{i!(n-i+1)!}+\frac{nB_{n+1}}{n+1}.
$$

\bigskip

\section{Формула Туэнтера и теоремы взаимности}

\bigskip

\begin{prop}\label{prop2}
$$
g(a,b)\sum_{k=0}^{b-1}\frac{ka}{b}h\left(\frac{ka}{b}\right)+g(b,a)\sum_{k=0}^{a-1}\frac{kb}{a}f\left(\frac{kb}{a}\right)
$$
$$
=g(a,b)\sum_{k=0}^{b-1}\frac{k}{b}h\left(\frac{k}{b}\right)+g(b,a)\sum_{k=0}^{a-1}\frac{k}{a}f\left(\frac{k}{a}\right)
$$
$$
+(g(a,b)\sum_{k\in C_{a,b}}\left(\left(\frac{k}{b}+1\right)h\left(\frac{k}{b}+1\right)-\frac{k}{b}h\left(\frac{k}{b}\right)\right)
$$
$$
+g(b,a)\sum_{k\in C_{a,b}}\left(\left(\frac{k}{a}+1\right)h\left(\frac{k}{a}+1\right)-\frac{k}{a}h\left(\frac{k}{a}\right)\right),
$$
в частности, для $1$-периодической функции $h$ имеем
$$
g(a,b)\sum_{k=0}^{b-1}\frac{ka}{b}h\left(\frac{ka}{b}\right)+g(b,a)\sum_{k=0}^{a-1}\frac{kb}{a}f\left(\frac{kb}{a}\right)
$$
$$
=g(a,b)\sum_{k=0}^{b-1}\frac{k}{b}h\left(\frac{k}{b}\right)+g(b,a)\sum_{k=0}^{a-1}\frac{k}{a}f\left(\frac{k}{a}\right)
$$
$$
+\sum_{k\in C_{a,b}}\left(g(a,b)h\left(\frac{k}{b}\right)+g(b,a)h\left(\frac{k}{a}\right)\right),
$$
для $g(a,b)=b^n$, $h(x)=\bar B_n(x)$ имеем
$$
ab^ns_n(a,b)+ba^ns_n(b,a)=\sum_{i=0}^n\binom{n}{i}(a^i+b^i)B_iS_{a,b}(n-i)
$$
$$
+\sum_{i=0}^{n+1}\frac{n!(a^{n-i+1}+b^{n-i+1})B_iB_{n-i+1}}{i!(n-i+1)!}+\frac{2nB_{n+1}}{n+1}.
$$
\end{prop}

  Доказательство. Применяем к суммам формулу Туэнтера (5) и складываем результаты с коэффициентами $g(a,b)$ и $g(b,a)$. Для $g(a,b)=b^n$, $h(x)=\bar B_n(x)$, используя равенства (2) и (11), получим
$$
ab^ns_n(a,b)+ba^ns_n(b,a)=b^n\sum_{k=1}^{b-1}\frac{k}{b}B_n\left(\frac{k}{b}\right)+a^n\sum_{k=1}^{a-1}\frac{k}{a}B_n\left(\frac{k}{a}\right)
$$
$$
+\sum_{k\in C_{a,b}}\left(b^nB_n\left(\frac{k}{b}\right)+a^nB_n\left(\frac{k}{a}\right)\right)
$$
$$
=\sum_{i=0}^{n+1}\frac{n!(b^{n-i+1}+a^{n-i+1})B_iB_{n-i+1}}{i!(n-i+1)!}+\frac{2nB_{n+1}}{n+1}
$$
$$
+\sum_{k\in C_{a,b}}\sum_{i=0}^n\binom{n}{i}(b^i+a^i)B_ik^{n-i}.
$$

\bigskip

\section{Степенные суммы спектральных чисел}

\bigskip

Мы получим две формулы для (сдвинутых) степенных сумм спектральных чисел особенности Брискорна $z_1^a+z_2^b$. Для первой формулы мы приведём краткое доказательство, но не будем выписывать это громоздкое равенство полностью. Более компактная формула представлена в Предложении 5. Нам понадобится известное равенство для степенных сумм
$$
\sum_{k=0}^{m-1}(x+k)^n=\frac{B_{n+1}(x+m)-B_{n+1}(x)}{n+1}.            \eqno (12)
$$
Применяя формулу (10), получим
$$
(ab)^n\sum_{i=1}^{b-1}\sum_{j=1}^{a-1}\left(\frac{x}{ab}+\frac{j}{a}+\frac{i}{b}\right)^n
$$
$$
=\sum_{\substack{0<i<b,0<j<a\\ia+jb<ab}}(x+ia+jb)^n+\sum_{\substack{0<i<b,0<j<a\\ia+jb<ab}}(x+2ab-ia-jb)^n.                       \eqno (13)
$$
Согласно \cite{6}
$$
\sum_{\substack{0\le i<b,0\le j<a\\ia+jb<ab}}(x+ia+jb)^{n-2}=\frac{B_{n-1}(x+ab)}{n-1}
$$
$$
-\frac{1}{n(n-1)ab}\sum_{k=0}^n\binom{n}{k}a^{n-k}b^k(B_{n-k}(b)-B_{n-k})B_k\left(\frac{x}{b}\right).                \eqno (14)
$$
Суммы в правой части формулы (13) отличаются от суммы вида (14) на суммы вида (12).

\begin{prop}\label{prop2}
$$
\sum_{i=1}^{b-1}\sum_{j=1}^{a-1}\left(\frac{x}{ab}+\frac{j}{a}+\frac{i}{b}\right)^n=\frac{x^n}{(ab)^n}
$$
$$
-\frac{B_{n+1}\left(\frac{x}{a}+b\right)-B_{n+1}\left(\frac{x}{a}\right)}{b^n(n+1)}-\frac{B_{n+1}\left(\frac{x}{b}+a\right)-B_{n+1}\left(\frac{x}{b}\right)}{a^n(n+1)}
$$
$$
+\frac{B_{n+1}(x+ab)-B_{n+1}(x)}{(ab)^n(n+1)}+\frac{1}{(ab)^n}\sum_{i=0}^n\binom{n}{i}((x+ab)^{n-i}-x^{n-i})S_{a,b}(i),
$$
в частности, для $x=0$ 
$$
\sum_{i=1}^{b-1}\sum_{j=1}^{a-1}\left(\frac{j}{a}+\frac{i}{b}\right)^n=\sum_{i=0}^n\binom{n}{i}\frac{S_{a,b}(i)}{(ab)^i}
$$
$$
+\frac{1}{n+1}\left(\frac{B_{n+1}(ab)}{(ab)^n}-\frac{B_{n+1}(a)}{a^n}-\frac{B_{n+1}(b)}{b^n}-B_{n+1}\left(\frac{1}{(ab)^n}-\frac{1}{a^n}-\frac{1}{b^n}\right)\right).
$$
\end{prop}

  Доказательство. Используя формулу Морделла (8) и равенство (12), получим
$$
(ab)^n\sum_{i=1}^{b-1}\sum_{j=1}^{a-1}\left(\frac{x}{ab}+\frac{j}{a}+\frac{i}{b}\right)^n=x^n-a^n\sum_{i=0}^{b-1}\left(\frac{x}{a}+i\right)^n-b^n\sum_{j=0}^{a-1}\left(\frac{x}{b}+j\right)^n
$$
$$
+\sum_{k=0}^{ab-1}(x+k)^n+\sum_{k\in C_{a,b}}((x+ab+k)^n-(x+k)^n)
$$
$$
=x^n-a^n\frac{B_{n+1}\left(\frac{x}{a}+b\right)-B_{n+1}\left(\frac{x}{a}\right)}{n+1}-b^n\frac{B_{n+1}\left(\frac{x}{b}+a\right)-B_{n+1}\left(\frac{x}{b}\right)}{n+1}
$$
$$
+\frac{B_{n+1}(x+ab)-B_{n+1}(x)}{n+1}+\sum_{k\in C_{a,b}}\sum_{i=0}^n\binom{n}{i}((x+ab)^{n-i}-x^{n-i})k^i.
$$

\bigskip

\section{Интеграл от сигнатуры торического узла}

\bigskip

  В \cite{3} для $(a,b)$-торического узла доказано равенство 
$$
\int_0^1\sigma(e^{2\pi ix})dx=\sum_{i=1}^{b-1}\sum_{j=1}^{a-1}\left(2\left|\frac{j}{a}+\frac{i}{b}-1\right|-1\right)        \eqno (15)
$$
Мы хотим вывести из этого равенства формулу (9). Используя тот факт, что число слагаемых в приведённых ниже суммах равно $(a-1)(b-1)/2$ \cite{14} и формулу (10), получим следующее продолжение равенства (15)
$$
(15)=\sum_{\substack{0<i<b,0<j<a\\ia+jb<ab}}\left(2\left|\frac{j}{a}+\frac{i}{b}-1\right|-1\right)+\sum_{\substack{0<i<b,0<j<a\\ia+jb<ab}}\left(2\left|2-\frac{j}{a}-\frac{i}{b}-1\right|-1\right)
$$
$$
=\sum_{\substack{0<i<b,0<j<a\\ia+jb<ab}}\left(1-2\left(\frac{j}{a}+\frac{i}{b}\right)\right)+\sum_{\substack{0<i<b,0<j<a\\ia+jb<ab}}\left(1-2\left(\frac{j}{a}+\frac{i}{b}\right)\right)
$$
$$
=(a-1)(b-1)-4\sum_{\substack{0<i<b,0<j<a\\ia+jb<ab}}\left(\frac{j}{a}+\frac{i}{b}\right).
$$
Согласно \cite{14}
$$
\sum_{\substack{0<i<b,0<j<a\\ia+jb<ab}}\left(\frac{j}{a}+\frac{i}{b}\right)=\frac{(a-1)(b-1)}{3}+\frac{(a-1)(b-1)(a+b+1)}{12ab},
$$
что доказывает формулу (9).

\end {document}